\documentclass[12pt]{article}
\usepackage{amsfonts,amsmath}
\usepackage{amssymb,latexsym}
\usepackage[dvips]{epsfig}
\usepackage{amscd}
 
\title{Hopfological algebra and categorification at a root of 
unity: the first steps}
\author{ Mikhail Khovanov}
\date{September 6, 2005}
 
\newtheorem{prop}{Proposition}
\newtheorem{theorem}{Theorem}
\newtheorem{lemma}{Lemma}

\newtheorem{corollary}{Corollary}

\newcommand{\oplusop}[1]{{\mathop{\oplus}\limits_{#1}}}

\begin{document} 
\maketitle
\baselineskip 14pt

\def\C{\mathbb C}
\def\R{\mathbb R}
\def\N{\mathbb N}
\def\Z{\mathbb Z}
\def\Q{\mathbb Q}
\def\l{\lbrace}
\def\r{\rbrace}
\def\lra{\longrightarrow}

\newcommand{\define}{\stackrel{\mathrm{def}}{=}}
\def\und{\underline}
\newcommand{\Id}{\mathrm{Id}}

\newcommand{\umod}{-\underline{\mathrm{mod}}}
\newcommand{\ugmod}{-\underline{\mathrm{gmod}}}
\newcommand{\dmod}{\mathrm{-mod}}
\newcommand{\gmod}{\mathrm{-gmod}}

\def\tr{\mathrm{tr}}
\def\binom#1#2{\left( \begin{array}{c} #1 \\ #2 \end{array}\right)}
\def\Hom{\textrm{Hom}}
\def\mc{\mathcal} 
\def\mf{\mathfrak}
\def\vsp{\vspace{0.1in}}

\tableofcontents

\begin{abstract} 
Any finite-dimensional Hopf algebra H is Frobenius and the stable 
category of H-modules is triangulated monoidal. To H-comodule 
algebras we assign triangulated module-categories over the stable 
category of H-modules. These module-categories are generalizations
of homotopy and derived categories of modules over a 
differential graded algebra. We expect that, for suitable H, our 
construction could be a starting point in the program of categorifying
quantum invariants of 3-manifolds. 
\end{abstract}


\section{Hopfological algebra} \label{hopfological}

{\bf Introduction.} Two of the better known examples of triangulated 
categories are the category of complexes of modules over a 
ring $A,$ up to chain homotopy, and the derived category 
$D(A\dmod).$ The first category is the quotient of the abelian 
category of complexes of $A$-modules by the ideal of null-homotopic 
morphisms. The derived category is the localization of 
the homotopy category given by inverting quasi-isomorphisms (maps 
that induce isomorphisms on cohomology). 

Stable categories provide more exotic examples of 
triangulated categories. Recall that a finite-dimensional 
algebra $H$ over a field $k$ is called \emph{Frobenius} if it has 
a non-degenerate trace $\tr: H\to k.$ A (left) module 
over a Frobenius algebra is injective if and only if it is projective. 
The stable category $H\umod$ 
is the quotient of the category of $H$-modules by morphisms 
that factor through an injective module. This category is triangulated
(see [Ha]). A finite-dimensional Hopf algebra is Frobenius, 
and its stable category carries both triangular and monoidal 
structures. 

In this paper we consider a new class of triangulated 
categories, associated to comodule-algebras over finite-dimensional 
Hopf algebras. This class includes stable categories of modules over 
Hopf algebras, as well as homotopy and derived categories of 
modules over algebras. We refer to the study of 
triangulated categories from this class as \emph{Hopfological algebra,} for it 
being a mix of homological algebra and the theory of Hopf algebras. 

We were originally inspired by the problem of categorifying quantum 
invariants of 3-manifolds. In Sections~\ref{cyclotomic} and \ref{done} 
of this paper we suggest 
an approach to this problem via triangulated module-categories 
over suitable stable categories $H\umod.$ A topologically inclined reader
might prefer to 
jump to these sections after reading the next subsection.

\vspace{0.15in} 

{\bf Hopf algebra and its stable category.} 
Let $H$ be a finite-dimensional Hopf algebra over a field $k.$ 
A left integral $\Lambda\subset H$ is a nonzero element of $H$ 
with the property $h\Lambda = \epsilon(h) \Lambda$ for any $h\in H,$
where $\epsilon$ is the counit of $H.$ A theorem of Larson and 
Sweedler [LS], [Mo, Section 2] says that a left integral exists and is 
unique up to multiplication by invertible elements of the ground 
field $k.$ From now on we fix a left integral $\Lambda.$ 

The category $H\dmod$ of left $H$-modules 
is monoidal, with 
$H$ acting on the tensor product $M\otimes N$ of $H$-modules 
via the comultiplication $\Delta.$ 
Unless specified otherwise, all tensor products 
are taken over $k.$

\begin{prop} $H$ is a Frobenius algebra. An $H$-module is projective 
if and only if it is injective. 
\end{prop} 

Proof of the first claim can be found in [LS] or [Mo, Section 2]. 
For a proof that projective modules over a Frobenius algebra 
coincide with injective  modules see [L, Chapters 3, 15, 16] 
(most of the other textbooks that cover Frobenius algebras prove 
this result only for finite-dimensional modules). 

Define the stable category $H\umod$ as the quotient of $H\dmod$ 
by the ideal of morphisms 
that factor through an injective (=projective) module. 
This category is triangulated, with the shift functor given by $TM= I/M$ 
where $M\subset I$ is an inclusion of $M$ as a submodule of 
an injective module. 
In general, for any Frobenius algebra $B,$ the stable category 
$B\umod$ has a triangular structure, with the shift functor as above, see [Ha]. 
The stable category $B\umod$ is trivial if and only if $B$ is semisimple.  

\begin{prop} \label{prop-tens} 
The tensor product $H\otimes M$ of the free $H$-module $H$ 
and an $H$-module $M$ is free  of rank equal the 
dimension of $M.$ The tensor product $P\otimes M$ of a 
projective $H$-module $P$ with an $H$-module $M$ is projective. 
\end{prop} 

The first statement is very well-known, see [Mo] or any other 
textbook on Hopf algebras, or proof of lemma~\ref{threemaps} 
below. The second statement follows at once from the first. 
$\square$  

In particular, any $H$-module $M$ admits a canonical embedding 
into a projective module, via the map $M\to H\otimes M$ given 
by  $m\mapsto \Lambda\otimes m.$ Notice that $H\Lambda=k\Lambda$ 
is a one-dimensional submodule of $H.$ Thus, we can define the 
shift functor 
$$H\umod\hspace{0.1in} \stackrel{T}{\lra} \hspace{0.1in} H\umod $$ 
by $TM= (H/k\Lambda)\otimes M.$ 

Proposition~\ref{prop-tens} also implies that tensor product 
of $H$-modules descends to a bifunctor in the stable 
category $H\umod$ and makes the stable category monoidal. 
The triangular and monoidal structures of $H\umod$ are 
compatible. For instance, tensoring with an $H$-module takes 
distinguished triangles to distinguished triangles, and 
there are isomorphisms
 $$ (TM) \otimes N \cong T (M\otimes N) \cong M \otimes TN.$$
We will say, informally, that $H\umod$ is a triangulated 
monoidal category, without giving a definition of the latter.

The category $H\umod$ has a symmetric monoidal structure for 
cocommutative $H.$ 
Only this case seems to have been thoroughly studied. 
Stable categories $k[G]\umod$ 
appear in the modular representation theory of finite groups $G,$ 
see [Ca], [BCR]. Voevodsky's working definition of a 
symmetric monoidal triangulated category [V], [MVW] was refined 
by May [Ma] (also see [KN]). 
For a relation between symmetric monoidal triangulated categories
and algebraic topology we refer the reader to Strickland [St] 
and references therein.  

\vspace{0.15in} 

{\bf Comodule-algebras and module-categories.} 
A left $H$-comodule algebra $A$ is an associative 
$k$-algebra equipped with a map
$$\Delta_A : A \lra H \otimes A $$ 
making $A$ an $H$-comodule and such that $\Delta_A$ is a 
map of algebras (where $H\otimes A$ has the product algebra 
structure). Equivalently, the following identities hold:  
\begin{eqnarray*} & & 
(\epsilon\otimes \Id)\Delta_A = \Id, \hspace{0.2in} 
  (\Delta\otimes \Id)\Delta_A = (\Id\otimes \Delta_A)\Delta_A, \\
& & \Delta_A(1)=1\otimes 1, \hspace{0.3in} 
 \Delta_A(ab)=\Delta_A(a)\Delta_A(b).
\end{eqnarray*} 
Let $A\dmod$ be the category of left $A$-modules. The tensor 
product $V\otimes M$ of an $H$-module $V$ and an $A$-module $M$ 
is naturally an $A$-module, via $\Delta_A.$ Therefore, this 
tensor product gives rise to a bifunctor 
$$H\dmod \hspace{0.05in} \times \hspace{0.05in} A\dmod 
  \lra A\dmod$$
compatible with the monoidal structure of $H\dmod,$ 
and makes $A\dmod$ into a module-category over $H\dmod.$ 

Let  $A_H\umod$ be 
the quotient of $A\dmod$ by the ideal of morphisms that factor through 
an $A$-module of the form $H\otimes N.$ In more detail, we call 
a morphism $f: M_1\to M_2$ of $A$-modules \emph{null-homotopic} if there 
exists an $A$-module $N$ and a decomposition $f=g_2 g_1$ with 
$$ M_1 \stackrel{g_1}{\lra} H\otimes N \stackrel{g_2}{\lra} M_2. $$  
The set of null-homotopic morphisms forms an ideal of $A\dmod.$ The 
quotient category $A_H\umod$ by this ideal has the same objects as 
$A\dmod,$ while the $k$-vector space of morphisms in 
$A_H\umod$ between objects $M_1, M_2$ is the quotient of 
$\Hom_{A\dmod}(M_1,M_2)$ 
by the subspace of null-homotopic morphisms. 

\vspace{0.1in} 

It's clear that $A_H\umod$ is a module-category over $H\umod,$ for 
the tensor product descends to a bifunctor 
$$ H\umod\hspace{0.05in} \times \hspace{0.05in} 
  A_H\umod \lra A_H\umod$$ 
compatible with the monoidal structure of $H\umod.$ 

\vspace{0.1in} 

{\bf Examples:} 

(a) If $A=H,$ with the standard $H$-comodule structure (via 
the comultiplication in $H$), then $H_H\umod$ is equivalent to $H\umod.$ 

(b) $A\subset H$ is a Hopf subalgebra, and 
$\Delta_A: A \to A\otimes A \subset H\otimes A$ is the 
comultiplication in $A.$ Since $H$ is a free $A$-module, 
the categories $A_H\umod$ and  $A\umod$ are equivalent.  

(c) If $A$ is a trivial $H$-comodule algebra, meaning 
$\Delta_A(a)=1\otimes a$ for all $a\in A,$ then the category 
$A_H\umod$ is trivial (any object is isomorphic to the zero object). 

(d) If $A$ is a semisimple $k$-algebra, the category $A_H\umod$ is 
trivial. 

(e) If $H$ is semisimple then $A_H\umod$ is trivial for any $A,$
since $H\otimes M$ contains $M$ as a direct summand, for any 
$A$-module $M.$ The 
composition  
 $$ M \stackrel{\Lambda\otimes \Id}{\lra} H \otimes M 
\stackrel{\epsilon\otimes \Id}{\lra}  M$$
is a multiple of the identity (for any $H$), and not equal 
to zero exactly when $H$ is semisimple.  

(f) Let $B$ be any $k$-algebra and form $A=H\otimes B,$
the tensor product of $H$ and $B,$ 
with the comodule structure $\Delta_A(h\otimes b) = \Delta(h)\otimes b.$ 
Then $A$ is a comodule-algebra over $H,$ with, in general,
highly non-trivial module-category $A_H\umod.$

(g) Suppose $B$ is a right $H$-module algebra. This means 
$B$ is a right $H$-module, and the multiplication and unit maps 
of $B$ are $H$-module maps (equivalently, $B$ is a left 
$H^{\ast}$ comodule-algebra). The smash product algebra 
$A=H\# B$ is the $k$-vector space $H\otimes B$ with 
the multiplication 
$$(h\otimes b)(l\otimes c)=\sum hl_{(1)} \otimes
  (b\cdot l_{(2)})c$$ 
and the left $H$-comodule structure 
$\Delta_A(h\otimes b) =\Delta(h)\otimes b.$  
It's easy to see that $A$ is a left $H$-comodule algebra. 

Example (f) is a special case of (g) when the right $H$-module 
structure on $B$ is trivial (factors through $\epsilon$). 
Of most interest to us is the example (g). We refer the reader 
to [Mo] or [DNR] for a generalization 
of (g) to crossed products and other extensions.  

\vspace{0.15in} 

{\bf The shift functor.} 
We next define the \emph{shift} or \emph{translation} endofunctor $T$ on 
the category $A_H\umod.$ Form the left $H$-module $H/(k\Lambda)$ 
and define, for any $M,$ 
$$T(M) = (H/(k\Lambda)) \otimes M.$$ 
 $T$ is a functor in the category of $A$-modules and in the category 
$A_H\umod$ as well.  Next, define the functor $T'$ in $A_H\umod$ by 
$$T'(M) =\mathrm{ker}(\epsilon) \otimes M.$$ 
The submodule $\mathrm{ker}(\epsilon)\subset H$ has codimension $1.$ 

\begin{prop} $T$ and $T'$ are mutually-inverse invertible functors 
in \newline $A_H\umod.$ \end{prop} 

\emph{Proof:} We have $T'T(M) = \mathrm{ker}(\epsilon) \otimes 
 (H/(k\Lambda)) \otimes M.$
To show that $T'T\cong \Id$ it suffices to decompose
the $H$-module $\mathrm{ker}(\epsilon) \otimes (H/H\Lambda )$ 
into a direct sum of a projective $H$-module $P$ and 
the trivial one-dimensional $H$-module $\underline{k}$ (on which 
$H$ acts via $\epsilon$), for then 
$$T'T(M) \cong (P\otimes M) \oplus(\underline{k}\otimes M)\cong 
 (P\otimes M) \oplus M,$$ 
and $P\otimes M$ is isomorphic to the zero object in the 
category $A_H\umod.$ 

The module $H\otimes (H/k\Lambda)$ is isomorphic to the free 
module $H^{n-1},$ where $n$ is the dimension of $H.$ Consider 
the diagram of surjective maps 
$$ H \stackrel{g}{\lra} H/k\Lambda  \stackrel{\epsilon\otimes \Id}
{\longleftarrow}  H\otimes (H/k\Lambda) \cong H^{n-1}$$
with $g$ being the obvious quotient map.  
Since $H$ is a free $H$-module, there exists a map $h: H \to H^{n-1}$ 
lifting $g,$ that is, $g= (\epsilon\otimes \Id)h.$ 

We claim that $h$ is injective, for otherwise it descends to
$$g': H/k\Lambda \lra H^{n-1}$$ 
such that $f g' = \Id.$ This would make $H/k\Lambda$ 
a direct summand of $H^{n-1},$ and a projective module. 
Then the short exact sequence 
$$0 \lra \underline{k} \stackrel{\Lambda}{\lra} H \lra H/k\Lambda 
 \lra 0$$ 
would split and the trivial one-dimensional $H$-module $\underline{k}$ 
would be projective and $H$ semisimple, the trivial case. 

Thus, $h: H \lra H^{n-1}$ is injective. Since $H$ is injective, we 
can find a splitting $H^{n-1}\cong h(H) \oplus P$ such that 
$f(P)=0.$ Therefore, 
$$\mathrm{ker}(\epsilon) \otimes (H/H\Lambda)\cong\underline{k}\oplus P$$
and $T'T\cong \Id.$ Likewise, $TT'\cong \Id.$    $\square$

\begin{lemma} \label{factor} 
 A null-homotopic $A$-module map $f:M_1\lra M_2$ 
  factors through the map $M_1 \stackrel{\Lambda\otimes \Id}{\lra}
  H\otimes M_1.$ 
\end{lemma} 

\emph{Proof:} It suffices to consider the case $M_2\cong H\otimes N.$ 
Since $H\otimes H\cong H^n$ as left $H$-modules and $H$ is injective, 
the inclusion 
$$\Lambda\otimes \Id : H \lra H\otimes H$$
admits a section $g: H\otimes H \lra H,$ 
an $H$-module map so that $g(\Lambda\otimes \Id)=\Id.$ The 
composition 
$$M_1\stackrel{\Lambda\otimes \Id}{\lra} H\otimes M_1 
\stackrel{\Id \otimes f}{\lra} H\otimes H\otimes N 
\stackrel{g\otimes \Id}{\lra} H\otimes N$$  
equals $f.$ $\square$ 

\vspace{0.1in} 

\emph{Remark:} David Radford pointed out an explicit formula for $g.$ 
Choose any $k$-linear map $t: H \lra H$ such that $t(\Lambda)=1.$ 
Then 
$$g(a\otimes b)= \sum b_{(2)} t (s^{-1}(b_{(1)})a)$$ 
is a map $H\otimes H\lra H$ of left $H$-modules such that 
$g(\Lambda\otimes b)=b$ for all $b\in H.$ Here $s$ is the antipode 
of $H.$  

\vspace{0.1in} 

{\bf Triangles.} 
For an $A$-module morphism $u:X\to Y$ denote by 
$\underline{u}$ its residue in the quotient category $A_H\umod.$ 
The cone $C_u$ of $u$ is defined as the pushout 
 \begin{equation*}
  \begin{CD}
     X @>{u}>> Y    \\
     @V{x}VV      @V{v}VV \\
     H\otimes X @>{\overline{u}}>> C_u   \end{CD}
 \end{equation*}
where $x=\Lambda\otimes \Id,$ $\overline{x}$ (see below) is the 
quotient map, so that we have a diagram of short exact 
sequences in $A\dmod$
 \begin{equation*} 
   \begin{CD}
   0 @>>> X @>{x}>> H\otimes X @>{\overline{x}}>> TX @>>> 0 \\
   @VVV    @V{u}VV    @V{\overline{u}}VV    @V{=}VV   @VVV  \\
   0 @>>> Y @>{v}>> C_u @>{w}>> TX @>>> 0
   \end{CD}
 \end{equation*} 
In this subsection we use notations $\overline{u}, \underline{u},$ 
etc. from Happel [Ha].

We define a \emph{standard distinguished triangle} as a sextuple 
$$ X \stackrel{u}{\lra} Y \stackrel{v}{\lra} C_u \stackrel{w}{\lra}TX$$ 
associated with a morphism $u.$ If $u$ and $u'$ are homotopic 
(i.e. $\underline{u}=\underline{u}'$), sextuples 
assigned to $u$ and $u'$ are isomorphic as diagrams of 
objects and morphisms in $A_H\umod.$ Thus, to a morphism 
$\underline{u}$ in $A_H\umod$ we can assign a sextuple 
$$  X \stackrel{\underline{u}}{\lra} Y \stackrel{\underline{v}}{\lra} 
C_u \stackrel{\underline{w}}{\lra}TX$$ 
in $A_H\umod.$ A sextuple 
$ X \stackrel{\underline{u}}{\lra} Y \stackrel{\underline{v}}{\lra} 
 Z \stackrel{\underline{w}}{\lra}TX$ of objects and morphisms in 
$A_H\umod$ is called a \emph{distinguished triangle} if it is isomorphic 
in $A_H\umod$ to a standard distinguished triangle.  

\begin{theorem} The category $A_H\umod$ is triangulated, with 
the shift functor $T$ and distinguished triangles as defined above. 
\end{theorem} 

The proof is almost identical to that of Theorem 2.6 in [Ha]. 
We need to check that the axioms TR1-TR4 of a triangulated 
category hold for $A_H\umod.$ The axiom TR1 is obvious. Let's 
prove TR3 following Happel [Ha], with minor modifications and 
using his notations whenever possible. The axiom TR3 
says that, given two triangles $(X,Y,Z,u,v,w)$ and $(X',Y',Z',u', 
 v', w')$ are morphisms $f:X\to X',$ $g:Y\to Y'$ such that 
$u'f=gu$ in $A_H\umod,$ there exists a morphism $(f,g,h)$ from the first 
triangle to the second.

Clearly, it's enough to verify this axiom for standard triangles. 
We have two diagrams 
 \begin{equation*} 
   \begin{CD}
   X @>{x}>> H\otimes X @>{\overline{x}}>> TX   \\
   @V{u}VV    @V{\overline{u}}VV    @V{=}VV     \\
   Y @>{v}>> C_u @>{w}>> TX
   \end{CD}
 \end{equation*} 
 \begin{equation*} 
   \begin{CD}
   X' @>{x'}>> H\otimes X' @>{\overline{x}'}>> TX'   \\
   @V{u'}VV    @V{\overline{u}'}VV    @V{=}VV     \\
   Y' @>{v'}>> C_{u'} @>{w'}>> TX'
   \end{CD}
 \end{equation*} 
and morphisms $f,g$ such that $\underline{u'f}=\underline{gu}.$ 
By lemma~\ref{factor}, there exists 
$$\alpha: H\otimes X \lra Y'$$ 
such that $gu=u'f+\alpha x.$ The map 
$$I_f= \Id \otimes f : H\otimes X \to H\otimes X'$$ 
satisfies $x'f=I_f x.$ For the map $Tf: TX\to TX'$ 
induced by $f$ we have $\overline{x}' I_f = Tf\hspace{0.03in} 
\overline{x}.$ 
The morphisms $v'g: Y \to C_{u'}$ and $\overline{u'} I_f + 
 v'\alpha: H\otimes X \to C_{u'}$ satisfy 
$$v' gu = v'u'f + v'\alpha x= \overline{u}' x'f+ v'\alpha x= 
(v'\alpha + \overline{u}' I_f) x.$$
Since $C_u$ is a pushout, there exists a morphism $h:C_u\to C_{u'}$ 
such that 
$$hv=v' g, \hspace{0.15in} h \overline{u} = \overline{u}' I_f + 
 v'\alpha.$$ 
Furthermore, $w' h v = Tf\hspace{0.03in} wv,$ since 
$Tf\hspace{0.03in} wv=0$ and $w'h v = w' v' g =0.$ Moreover,  
$$  Tf\hspace{0.03in} w \overline{u}= Tf \hspace{0.03in} 
\overline{x} = \overline{x}' I_f= 
 w \overline{u}' I_f = w'\overline{u}' I_f + w' v'\alpha= 
 w'(\overline{u}' I_f+ v'\alpha) = w' h \overline{u}.$$ 
This implies $w' h = Tf\hspace{0.03in} w$ and 
$(\underline{f}, \underline{g}, 
 \underline{h})$ is a morphism of triangles. 

\vspace{0.1in} 

Axiom TR2 says that $(Y,Z,TX,v,w,-Tu)$ is a triangle 
if $(X,Y,Z,u,v,w)$ is. 
Happel's proof of this axiom transfers to our case without any 
changes (as long as we substitute $H\otimes X$ and $H\otimes Y$ 
for his $I(X)$ and $I(Y)$). 

\vspace{0.1in} 

The last axiom TR4 is the octahedral axiom. Happel's argument  
works for us as well (after substituting $H\otimes N$ for 
his $I(N)$ throughout for various $N$), once we observe that 
the sequence of $A$-modules 
$$Y \stackrel{i}{\lra} H \otimes C_u \lra 
(H\otimes C_u)/i(Y)$$ 
(where $i=\Lambda\otimes i_1$ and $i_1$ is the inclusion 
$Y\to C_u$) is isomorphic to the direct sum of 
$$Y \stackrel{y}{\lra} H\otimes Y \stackrel{\overline{y}}{\lra} TY$$ 
and 
$$ 0 \lra H\otimes (H/k\Lambda) \otimes X 
 \stackrel{\Id}{\lra}  H\otimes (H/k\Lambda) \otimes X.$$

$\square$ 

\vspace{0.1in} 

For any $H$-module $V,$  the inclusions 
$$ \Id\otimes \Lambda: V \lra V\otimes H\hspace{0.2in}
\mathrm{and} \hspace{0.2in}\Lambda\otimes\Id: V\lra H\otimes V $$ 
are $H$-module maps. 

\begin{lemma} \label{threemaps} 
There exists a functorial in $V$ isomorphism 
$r: V\otimes H\lra H\otimes V$ intertwining the above inclusions: 
 $r\cdot(\Id \otimes \Lambda) = \Lambda \otimes \Id.$
\end{lemma}

\emph{Proof:}  Let $r=r_3 r_2 r_1,$ where 
$$ V\otimes H \stackrel{r_1}{\lra} V'\otimes H \stackrel{r_2}{\lra} H\otimes V' 
 \stackrel{r_3}{\lra} H\otimes V$$ 
and 
\begin{eqnarray*} 
 r_1 (v\otimes h) & = & \sum s^{-1}(h_{(1)})v \otimes h_{(2)},     \\
 r_2 (v\otimes h) & = &   h g \otimes v, \\
 r_3 (h\otimes v) & = &  \sum  h_{(1)} \otimes h_{(2)} v .   
\end{eqnarray*} 
Here $V'=V$ as a $k$-vector space, but $H$ acts trivially on $V'.$ 
We denoted by  $g$  the distinguished grouplike element of $H,$ determined by the 
condition $x\lambda = x(g) \lambda$ for all $x\in H^{\ast},$ where 
$\lambda$ is a right integral for $H^{\ast},$ i.e. $\lambda x = x(1)\lambda$ 
(see [Ra] for details and [Ku] for a diagrammatical approach).   
The maps $r_1,r_2,r_3$ are functorial in V invertible $H$-module maps.  
The equation $r(\Lambda\otimes \Id) = \Id \otimes \Lambda$
reduces, by a straightforward simplification, to the equation 
$$ \sum \Lambda_{(1)} \otimes \Lambda_{(2)} = \sum s^2(\Lambda_{(2)}) g \otimes 
   \Lambda_{(1)}$$ 
which is part (d) of Theorem 3 in [Ra]. $\square$ 

\vspace{0.15in} 

An additive functor $F:\mc{C}\lra \mc{D}$ between triangulated categories
is called \emph{exact} if it takes distinguished triangles to 
distinguished triangles. The lemma implies that the tensor product 
$$ V\otimes X \lra V\otimes Y \lra V\otimes Z \lra V\otimes TX$$ 
of a distinguished triangle in $A_H\umod$ and an $H$-module $V$ 
is a distinguished triangle, so that tensoring with an $H$-module $V$ is 
an exact functor in $A_H\umod.$ We say, informally, that 
$A_H\umod$ is a triangulated module-category over $H\umod.$

The Grothendieck group $K(\mc{C})$ of a triangulated category 
$\mc{C}$ is the group with generators $[X],$ over all objects 
$X$ of $\mc{C},$ and relations 
$[Y]=[X]+[Z]$ whenever there exists a distinguished triangle 
$X\to Y\to Z\to TX.$ These relations imply $[TX]=-X.$ 

To have nontrivial Grothendieck groups, we should avoid 
"infinities" (for instance, allowing infinite direct sums leads 
to $[X]=0$ for any object $X$) and restrict to modules from 
a particular small class. Whenever we are considering 
Grothendieck groups, $H\umod$ will denote the stable category 
of finite-dimensional modules, while for objects 
of $A_H\umod$ we can take all finite-dimensional $A$-modules, 
or finitely-generated $A$-modules (if $A$ is Noetherian), or, 
more generally, any class of $A$-modules closed under 
taking finite direct sums, subquotients, and 
tensoring with $H.$ 

\vspace{0.15in} 

Since $H\umod$ has an exact tensor product, $K(H\umod)$ is 
a ring. As an abelian group, $K(H\umod)$ is finitely-generated 
and may have torsion. For quasi-triangular $H,$ the ring 
$K(H\umod)$ is commutative. 

\begin{corollary} 
The Grothendieck group $K(A_H\umod)$ is a left module over the 
ring $K(H\umod).$ 
\end{corollary} 

Suppose now that $A=H\#B$ is a smash product of $H$ and a right $H$-comodule 
algebra $B,$ see earlier example (g).  Then we have the restriction functor 
from $A$-modules to $H$-modules (due to $H$ being a subalgebra 
of $A$) which descends to an exact 
functor from $A_H\umod$ to $H\umod.$ We say that a morphism $f:M\to N$ 
in $A_H\umod$ is a \emph{quasi-isomorphism} if it restricts to an isomorphism in 
$H\umod.$ 

\begin{prop}\label{prop-loc}
 Quasi-isomorphisms in $A_H\umod$ constitute a localizing class. 
\end{prop} 

\emph{Proof:} The composition of quasi-isomorphisms is a quasi-isomorphism. 
Suppose $X\stackrel{s}{\lra} Y \stackrel{f}{\longleftarrow} Z$ is a
 diagram of objects and morphisms in $A\dmod$ with $s$ a quasi-isomorphism. 
Form the pull-back diagram
 
 \begin{equation*}
   \begin{CD}
     C @>{g}>> H\otimes Y  \\
     @V{h}VV      @VV{\epsilon\otimes \Id}V \\
     X\oplus Z @>{(s,f)}>> Y
   \end{CD}
 \end{equation*}
 
The maps $h_X: C \lra X$ and $h_Z: C\lra Z$ (where $h=(h_X,h_Z)$)
make the diagram below commutative in $A_H\umod$ 

\begin{equation} \label{a-cd} 
   \begin{CD}
     C @>{h_Z}>> Z  \\
     @V{-h_X}VV      @VV{f}V \\
     X @>{s}>> Y
   \end{CD}
 \end{equation}

\begin{lemma} The map $h_Z$ is a quasi-isomorphism. 
\end{lemma} 

\emph{Proof:} $h_Z$ is surjective as a map of $A$-modules and 
we have a short exact sequence 
$$ 0 \lra \mathrm{ker}(h_Z) \lra C \stackrel{h_Z}{\lra} Z \lra 0.$$
The module $\mathrm{ker}(h_Z)$ consists of pairs $(x,a)$ where 
$x\in X,$ $a\in H\otimes Y$ and $s(x)= (\epsilon \otimes \Id) a.$ 
Thus, this module is the kernel of the surjective quasi-isomorphism 
$$ X\oplus (H \otimes Y) \stackrel{(s,-\epsilon \otimes \Id)}{\lra} 
 Y.$$ 
But the kernel of a surjective quasi-isomorphism is a projective 
(hence, injective) $H$-module, and the above sequence splits, 
implying the lemma's statement. $\square$ 

\vspace{0.08in} 

We see that a diagram 
 $X\stackrel{s}{\lra} Y \stackrel{f}{\longleftarrow} Z$
with $s$ a quasi-isomorphism lifts to a commutative square 
(\ref{a-cd}) with $h_Z$ a quasi-isomorphism. 

\vspace{0.08in} 

A dual argument implies that any diagram 
 $X\stackrel{s}{\longleftarrow} Y \stackrel{f}{\longrightarrow} Z$
with $s$ a quasi-isomorphism extends to a commutative square 
\begin{equation} \label{b-cd} 
   \begin{CD}
     Y @>{s}>> X  \\
     @V{f}VV      @VV{g}V \\
     Z @>{t}>> W
   \end{CD}
 \end{equation}
with $t$ a quasi-isomorphism. 

\vspace{0.08in} 

Next, we need to show that, given a diagram 
$X\stackrel{f}{\lra} Y \stackrel{s}{\lra} Z$ in 
$A_H\umod$ with $s$ a quasi-isomorphism and $sf=0$ 
in $A_H\umod,$ there exists a quasi-isomorphism 
$W\stackrel{t}{\lra} X$ such that $ft=0.$ 

By 
adding to $Y$ a direct summand of the form $H\otimes U,$ 
if necessary, we can assume $sf=0$ in the abelian 
category $A\dmod$ rather than just in the quotient 
category $A_H\umod.$ Let 
\begin{eqnarray*}
 G_1 & = & X \oplus (H\otimes Y) \oplus (H\otimes H \otimes Z), \\
 G_2 & = & Y \oplus (H\otimes Z),
\end{eqnarray*} 
and $g:G_1 \to G_2$ be the sum of the following four maps:  
 \begin{eqnarray*} 
   s & :  & X\lra Y, \\
   \epsilon\otimes \Id & :  & H\otimes Y \lra Y, \\
    \Id \otimes s & : & H\otimes Y \lra H\otimes Z, \\
   \epsilon \otimes \Id^{\otimes 2} & : & H \otimes H \otimes Z 
   \lra H\otimes Z. 
 \end{eqnarray*} 
Let $W=\mathrm{ker}(g)$ and $W\stackrel{t}{\lra} X$ be 
the composition of inclusion and projection onto the first summand, 
$W \subset G_1 \rightarrow X.$ As a map of $k$-vector spaces, 
$t$ is surjective with the kernel a projective $H$-module (an easy 
exercise). Hence $t$ is a quasi-isomorphism. The composition 
$ft$ factors through $H\otimes Y;$ thus $ft=0$ in $A_H\umod.$ 

\vspace{0.1in} 

A dual argument produces a quasi-isomorphism 
$Y\stackrel{s}{\lra} Z$ such that $sf=0,$ given 
a diagram $W\stackrel{t}{\lra} X \stackrel{f}{\lra} Y$ 
with $t$ a quasi-isomorphism and $ft=0.$ Proposition~\ref{prop-loc} 
follows. 

$\square$ 

\vspace{0.15in} 

Recall that $A$ is the smash product $H\# B.$ 
Denote by $D(B,H)$ the localization of $A_H\umod$ with respect 
to quasi-isomorphisms. 
It has the same objects as $A_H\umod$ and morphisms from $M$ to $N$ 
are equivalence classes of diagrams 
$M \stackrel{s}{\longleftarrow} M' \stackrel{f}{\lra} N$ where $s$ is a
quasi-isomorphism. We refer the reader to [GM], [We] for properties 
 of localizations.

\begin{corollary} $D(B,H)$ is a triangulated category. 
Tensoring with an $H$-module is an exact functor in $D(B,H),$
and the Grothendieck group of $D(B,H)$ is a module over 
$K(H\umod).$  
\end{corollary} 

Throughout the paper we denote the cyclic group $\Z/n\Z$ by $\Z_n.$ 
If $H$ is a $\Z-\mathrm{graded}$ or $\Z_n-\mathrm{graded}$ Hopf 
algebra with a homogeneous left integral, the 
stable category  $H-\underline{\mathrm{gmod}}$
of graded $H$-modules is
triangulated monoidal. Its Grothendieck ring is a $\Z[q,q^{-1}]$-algebra, 
with the multiplication by $q$ corresponding to the grading shift. 
If $A$ is a graded $H$-comodule algebra, the category 
$A_H-\underline{\mathrm{gmod}}$ of graded $A$-modules modulo homotopies
is triangulated. 
All other results above generalize to the graded case as well. 

\vspace{0.1in} 

If $H$ is a finite-dimensional Hopf super-algebra (a Hopf algebra
in the category of super $k$-vector spaces), the stable category 
$H\umod$ is triangulated and all the arguments above  
work in this case as well (with $A$ being a comodule-superalgebra 
over $H,$ etc). Moreover, the construction generalizes 
to the situation when $H$ is a $\Z$-graded Hopf super-algebra, 
with the super $\Z_2$-grading given by the mod 2 reduction of 
the $\Z$-grading, and we consider the stable category of $\Z$-graded 
modules.


\section{Examples} 

{\bf DG algebras.} A differential graded $k$-algebra $B$ comes 
equipped with a degree 1 differential $D$ such that 
$D^2=0$ and 
$$D(ab)= D(a) b + (-1)^{|a|} a D(b), \hspace{0.3in} a,b\in B.$$
To see the link between DGAs and comodule-algebras, note that
the exterior algebra on one generator $H=k[D]/(D^2)$  
is a Hopf algebra in the category 
of $\Z$-graded super-vector spaces.

The smash product $A= B\# H$ is a right $H$-comodule algebra
via 
$$ A\stackrel{\Delta_A}{\lra} A\otimes H, \hspace{0.2in} 
  a\otimes h \longmapsto a \otimes \Delta(h). $$ 
As a vector space $A= B\otimes H,$ with the multiplication
$$ (b \otimes h)(b' \otimes h') = \sum b h_{(1)}(b')\otimes h_{(2)} h',$$ 
where $\Delta(h)= \sum h_{(1)} \otimes h_{(2)}.$ In particular, 
$\Delta(D) = D\otimes 1 + 1\otimes D,$ and
\begin{eqnarray*} 
 ( b\otimes 1)(b'\otimes h') & = & b b'\otimes h', \\ 
 ( b\otimes D)(b'\otimes h') & = & b D(b') \otimes h' + (-1)^{|b'|}
 b b' \otimes D h'. 
\end{eqnarray*} 
A differential graded left $B$-module is the same as a graded 
left $A$-module, so that the two abelian categories $B-\mathrm{dgmod}$ 
and $A\gmod$ (of DG modules 
over $B$ and graded $A$-modules, respectively) are equivalent. 
Next, the quotient category $A_H\umod$ is equivalent to the quotient 
of  $B-\mathrm{dgmod}$ by the ideal of null-homotopic morphisms. 
Finally, the derived category $D(B,H)$ is equivalent to the
derived category $D(B)$ of the dga $B,$ as defined, for instance, 
in [BL]. 
 
DG algebras constitute an indispensable tool in 
homological algebra. The above 
correspondence and results of the previous section suggest 
a generalization of the notion of a DG algebra and 
triangulated categories of modules assigned to it. 
A finite-dimensional Hopf algebra $H$ (more generally, 
a finite-dimensional graded Hopf (super) algebra) takes place of 
the homological equation $D^2=0.$ 
A module-algebra over $H$ takes place of a DG algebra.  
The homotopy category of differential graded $B$-modules generalizes to 
the category $A_H\umod,$ where $A$ is the smash 
product $H\#B,$ while $D(B,H)$ plays the role of the 
derived category of DG modules. We expect that large chunks of the 
theory of DG algebras will generalize to $H$-module algebras.

\vspace{0.15in} 

{\bf Finite groups.} Let $G$ be a finite group. $A$ is a
comodule-algebra over $k[G]$ if 
$$A=\oplusop{g\in G} A_g, \hspace{0.3in} A_gA_h\subset A_{gh}.$$
Then $\Delta_A(a) =g\otimes a$ for $a\in A_g.$ 
The category $A_{k[G]}\umod$ is a module-category over 
$k[G]\umod.$ We assume that $\mathrm{char}(k)$ 
divides the order of $G,$ for otherwise $A_{k[G]}\umod$ is trivial. 

If, in addition, there are elements $a_g\in A_g$ such that $a_1=1$ 
and $a_ga_h=a_{gh},$ then $k[G]\subset A,$ via the map 
$g\mapsto a_g,$ and $(\Delta_A)|_{k[G]}=\Delta,$ 
so we can form the derived 
category $D(A_1,k[G]).$ In this case $A_g=a_g A_1= A_1 a_g$ 
and $a_g x= \psi_g(x) a_g,$ for $x\in A_1$ and 
 an automorphism $\psi_g$ of the algebra $A_1.$ We see that 
$A$ is determined by $A_1$ and a homomorphism $\psi$ from $G$ 
to the automorphism group of $A_1.$ The algebra $A$ is then a 
smash product of $A_1$ and the group ring $k[G].$ The derived category 
$D(A_1,k[G])$ seems interesting already in the case of trivial 
$\psi,$ when $\psi_g=\Id$ for all $g$ and $A=k[G]\otimes A_1$ 
is the tensor product algebra. 

Although one rarely encounters $G$-graded rings for non-commutative 
groups $G,$ it's common for an algebra to be $\Z_n$-graded, for 
some $n.$ Let's specialize to $G=\Z_n$ and assume 
$k$ has finite characteristic which divides $n.$ 
The resulting categories $A_{k[\Z_n]}\umod$ and $D(A,k[\Z_n])$ 
resemble the categories we'll study in Section~\ref{cyclotomic} 
associated with the truncated polynomial Hopf algebra. 

The first nontrivial case $n=2$ and $\mbox{char }k=2$ produces 
categories similar to the usual homotopy and derived categories of 
modules over a ring. Let's pass to the additive notation for 
indices, so that $A=A_0\oplus A_1,$
and consider the simplest case $A=A_0\otimes k[\Z_2].$  
Denote by $g$ the generator of $\Z_2$ and let $d=g+1\in k[\Z_2].$ 
Then $d^2=0.$ A direct computation shows that an $A$-module $M$ is 
trivial in the category $A_{k[\Z_2]}\umod$ iff it's isomorphic, 
as an $A$-module, to $k[\Z_2]\otimes N,$ for some $A_0$-module $N.$ 
Likewise, an $A$-module $M$ is trivial in the derived category 
$D(A_0,k[\Z_2])$ iff $ker(d)=im(d).$ Thus, triangulated categories 
 $A_{k[\Z_2]}\umod$ and $D(A_0,k[\Z_2])$ are very close relatives 
of the homotopy category of complexes over the algebra $A_0$ and 
its derived category, respectively. Unlike the usual homotopy and 
derived categories, the modules here are not graded, and we have to 
work over a field of characteristic $2.$  

\vspace{0.15in} 

{\bf Bicomplexes.} Consider a 
$\Z\oplus\Z$-graded $\C$-algebra $\Lambda_2$ with generators 
$D_1, D_2$ of bidegrees $(1,0)$ and $(0,1),$ respectively, and 
relations 
$$D_1^2=0, \hspace{0.2in} D_2^2=0, \hspace{0.2in} D_1D_2+D_2D_1.$$ 
$\Lambda_2$ is simply the exterior algebra on two generators, and 
a graded Hopf algebra in the category of super-vector spaces, 
with $\Delta(D_i)=D_i \otimes 1 + 1\otimes D_i.$ Any $\Z\oplus 
\Z$-graded algebra $A$ with two super-commuting graded 
super-derivations $\partial_1, \partial_2$ (such that $\partial_i^2=0$) 
 is a $\Lambda_2$-comodule 
algebra, and one can form triangulated module-categories 
$A_{\Lambda_2}-\underline{\mathrm{gmod}}$ and $D_{A_{\Lambda_2}}(A)$ 
over $\Lambda_2-\underline{\mathrm{gmod}}.$ Let's call such $A$ a bi-DG-algebra.

Tensor products $A_1\otimes A_2,$ with each $A_i$ a DG-algebra,
deliver examples of bi-DG-algebras. 
If $X$ is a compact complex manifold, the Dolbeault bicomplex 
$$\Omega(X) = \oplusop{i,j}\hspace{0.05in}
\Omega^{i,j}(X), \hspace{0.1in} 
 \partial: \Omega^{i,j}(X) \lra \Omega^{i+1,j}(X), 
 \hspace{0.1in} \overline{\partial}: 
\Omega^{i,j}(X) \lra \Omega^{i,j+1}(X)$$ 
is a (super) commutative bi-DG-algebra. 

Denote by $\underline{\Omega(X)}$ the bicomplex $\Omega(X)$ 
viewed as an object of the stable category $\Lambda_2\ugmod.$ 
The Dolbeault and de 
Rham cohomology groups $H_{\partial}(X),$ $H_{\overline{\partial}}(X),$ 
$H_{\partial+\overline{\partial}}(X)$  can be recovered from 
$\underline{\Omega(X)},$ and are the image of $\underline{\Omega(X)}$ 
under suitable tensor functors from $\Lambda_2\ugmod$ to 
the abelian category of $\Z\oplus \Z$-graded (respectively, $\Z$-graded) 
vector spaces. 

For a general complex manifold $X$ the invariant  
$\underline{\Omega(X)}$ 
seems to carry more information than the cohomology groups 
and provides an example where 
the stable category of a Hopf algebra appears naturally. 
When $X$ is K\"ahler, the Hodge theory tells us that 
$\Omega(X)\cong P \oplus S$ where $P$ is an (infinite-dimensional) 
projective and $S$ a finite-dimensional semisimple 
$\Lambda_2$-modules, and 
$$S\cong \oplusop{i,j} H^{i,j}(X) \otimes S_{i,j},$$ 
with $S_{i,j}$ being the simple $\Lambda_2$-module which is 
$\C$ placed in bidegree $(i,j).$ In the stable category $P=0$ 
and $\underline{\Omega}(X) \cong S.$ Thus, for K\"ahler $X,$ 
the invariant $\underline{\Omega}(X)$ carries the same information 
as the Dolbeault cohomology.  

\vspace{0.2in} 

{\bf Finite type and tame Hopf algebras.} 
The stable category $H\umod$ serves as the target category 
for various exact functors from $A_H\umod$ and 
$D(B,H),$ where $A$ is a smash product, $A=H\# B.$ 
It seems much easier to investigate homological algebra 
over $H\umod$ when the objects of this base category  
could be classified. Both categories $H\dmod$ and $H\umod$ have 
the Krull-Shmidt property (we restrict now to finite-dimensional 
modules) and their objects split uniquely as direct sums 
of indecomposables. Since there are only finitely many 
indecomposable projective $H$-modules, $H$ has finite type iff 
$H\umod$ has finite type, and $H$ is a tame algebra iff $H\umod$ 
is tame. These two cases (when $H$ is finite or tame) seem 
particularly important, for then we can completely understand 
the objects of $H\umod.$  

A semisimple Hopf algebra $H$ has finite type, but then 
$H\umod$ is trivial. 
The following are examples of nonsemisimple Hopf algebras of finite type.  
\begin{itemize} 
\item The truncated polynomial algebra $k[X]/(X^n)$ where field 
$k$ has characteristic $p$ and $n=p^m,$ with $\Delta(X)=X\otimes 1 + 1 \otimes X$ 
(See Section~\ref{cyclotomic}). 
\item The restricted universal enveloping algebra of the p-Lie algebra with 
 a basis $\{x,y\}$ and relations $[y,x]=x,$ $x^{[p]}=0,$ $y^{[p]}=y.$ 
\item The group ring $k[\Z_n]$ of a cyclic group, where $p=\mathrm{char}(k)$
 divides $n.$ More generally, a group ring $k[G]$ has finite type if 
 $p$-Sylow subgroups of $G$ are cyclic.   
\item The Taft Hopf algebra (see Section~\ref{done}).  
\end{itemize} 
Examples of nonsemisimple Hopf algebras of tame type are  
\begin{itemize} 
\item the finite quantum group $U^{fin}_q(\mf{sl}_2)$ 
 when $q$ is a root of unity,   
\item the restricted 
universal enveloping algebra of the p-Lie algebra sl(2),
\item group rings $k[G],$ when $\mathrm{char}(k)=2$ and the 2-Sylow subgroup 
of $G$ is the Klein four group, 
\item cross-products of $\C[G]$ and the exterior algebra of $\C^2,$ 
for $G$ a finite subgroup of $SU(2),$
\item the Drinfeld double of the Taft Hopf algebra [EGST].  
\end{itemize} 
For a detailed study of finite and tame cocommutative Hopf 
algebras we refer the reader to [FS] and references therein. 

\vspace{0.1in} 

A common theme in modern algebraic topology and 
modular representation theory is to start with a complicated 
tensor triangulated category, such as the category of spectra 
or the stable category $k[G],$ for $G$ with a large p-Sylow 
subgroup, and then localize to a simpler tensor triangulated category
which is more accessible to computation [St]. 
We are suggesting here a different direction, where one starts 
with a rather uncomplicated tensor triangulated category, say $H\umod$ 
for a tame or finite type $H,$ but then goes on to study more 
sophisticated triangulated module-categories over this category. 

\vspace{0.1in} 


\section{A categorification of the cyclotomic ring} 
\label{cyclotomic} 

{\bf Category $\mc{H}'.$}
The Jones polynomial of links admits a categorification, 
being the Euler characteristic of a bigraded link homology 
theory [Kh]. The Jones polynomial also extends 
to a three-manifold invariant, the Witten-Reshetikhin-Turaev 
invariant [W], [RT], which has an interpretation via quantum sl(2) 
at a root of unity. For a closed three-manifold $M$ 
this invariant $\tau(M),$ which depends on an integer $n,$ 
takes values in the cyclotomic field $\Q(\zeta^{\frac{1}{4}}) 
 \subset \C,$
where $\zeta= \exp(\frac{2\pi i}{n}).$ More precisely, $\tau(M)$ 
lies in the ring $\Z[\zeta^{\frac{1}{4}},\frac{1}{n}]\subset 
\Q(\zeta^{\frac{1}{4}}).$ When $n$ is prime, suitably normalized 
$\tau(M)$ takes values in the ring of cyclotomic integers 
$\Z[\zeta]\subset \C,$ see [Mu], [MR], [Le]. 
Louis Crane and Igor Frenkel conjectured [CF] that some quantum 
invariants of 3-manifolds can be categorified. 

\vspace{0.1in} 

The Jones polynomial of a link in $\R^3$ lies in $\Z[q,q^{-1}].$ 
Categorification turns the Jones polynomial into a bigraded abelian 
group, and $q$ into a grading shift $\{1\}.$ When 
$q=\zeta$ is the primitive $n$-th root of unity, the naive 
approach to categorification would be to make the grading $n$-periodic, 
and work with $\Z_n$-graded objects. 

However, in the ring of cyclotomic integers (where quantum invariants 
tend to live) there is a stronger relation 
\begin{equation}\label{eq-sum} 
1+q+\dots + q^{n-1}=0,
\end{equation} 
 and that's the one we must categorify (at least when 
 $n$ is a prime, for composite $n$ a stronger relation 
 $\Psi_n(q)=0$ holds, where $\Psi_n$ is the $n$-th cyclotomic 
 polynomial). We let 
\begin{equation*} 
   R_n= \Z[q]/(1+q+\dots + q^{n-1})
\end{equation*} 
and try to categorify the ring $R_n.$ 

A categorification of a ring is an additive category $\mc{C}$
with an exact tensor product (this is slightly inaccurate, since to 
define exactness we need $\mc{C}$ to be either abelian or 
triangulated). An associative multiplication 
in the Grothendieck group $K(\mc{C})$ is given by 
 $$ [M][N] \define [M\otimes N] $$
where $[M]$ is the image of an object $M$ in the Grothendieck 
 group. The image $[{\bf 1}]$ of the unit object 
${\bf 1}$ of $\mc{C}$ is the unit element of the ring $K(\mc{C}).$   

For $K(\mc{C})$ to be commutative, it suffices to require 
that $\mc{C}$ is braided (has a compatible family of 
isomorphisms $M\otimes N \cong N\otimes M$). 
 For $K(\mc{C})$ to have a natural structure of a $\Z[q,q^{-1}]$-module, 
it suffices to have an exact invertible functor $\{ 1\}$ in $\mc{C},$ 
compatible with the tensor structure.  

To have a homomorphism $R_n \lra K(\mc{C})$ we need the relation  
 \begin{equation}\label{eq-grgr} 
  [{\bf 1} ] + [{\bf 1}\{ 1\}] + \dots + [{\bf 1}\{ n-1\}] =0.
 \end{equation} 
Positivity constraints seem to make this equation impossible 
in a abelian category. Equation (\ref{eq-grgr}) holds  
if $\mc{C}$ has an object $H,$ built out of $n$ shifts of 
the unit object ${\bf 1},$ such that $H$ is zero in the Grothendieck 
group, $[H]=0.$ 

A solution to this puzzle does exist in the realm of 
triangulated categories. Take the stable category $H\umod$ 
of finite-dimensional $H$-modules, where $H$ is a finite-dimensional 
Hopf algebra over a field $k.$ The stable category is triangulated 
and monoidal, and the Grothendieck group $K(H\umod)$ is a ring. 
The image of the $H$-module $H$ in the 
Grothendieck ring of the stable category is trivial, since 
$H$ is projective as a module over itself.  

We should avoid the case of semisimple $H,$ resulting in trivial 
stable category. 
One of the first examples of nonsemisimple finite-dimensional Hopf 
algebras is the truncated polynomial algebra 
$$ H = k[X]/(X^p),$$ 
where field $k$ has characteristic $p>0.$ The comultiplication and 
the antipode are given by $\Delta(X)=X\otimes 1 + 1\otimes X$ 
and $S(X)=-X.$ Furthermore, $H$ is graded, with $X$ in degree $1.$ 

Let $H\gmod$ be the category of finite-dimensional 
$\Z$-graded $H$-modules and $\mc{H}'=H\ugmod$ the corresponding stable 
category (of finite-dimensional graded $H$-modules modulo projectives). 
Both categories have symmetric tensor structure, 
and their Grothendieck groups are rings. 
The Grothendieck group $K(H\mathrm{-gmod})$ is 
isomorphic to $\Z[q,q^{-1}].$ 

\begin{prop}\label{cyclo}
 The Grothendieck group $K(\mc{H}')$ is naturally  
isomorphic to the ring of cyclotomic integers $R_p.$ Multiplication 
by $q$ corresponds to the grading shift. 
\end{prop} 

\emph{Proof:} Let $V_i = H/(X^{i+1} H),$ for $i=0,1, \dots, p-1.$
The module $V_{p-1} \cong H$ is free while $M_0$
is a simple $H$-module and the unit object of $H\mathrm{-gmod}$ 
and $\mc{H}'.$ We also denote $V_0$ by $\mathbf{1}.$ 
Existence of the Jordan form for  the nilpotent 
operator $X$ acting on an $H$-module implies that     
any module in $H\gmod$  decomposes uniquely (up to isomorphism) 
as a direct sum of modules $V_i\{ j\},$ for $i\le p-1$ 
and $j\in \Z.$ Here $\{j\}$ denotes the grading shift upward by $j.$
Likewise, any object of $\mc{H}'$ has a unique (up to 
isomorphism) decomposition as a direct sum of modules  $V_i\{ j\},$ for 
$i< p-1$ and $j\in \Z.$

In particular, any 
indecomposable object of $\mc{H}'$ is isomorphic to $V_i\{ j\}$ 
for a unique pair $(i,j)$ with $i<p-1.$ The module $V_i$ has a filtration 
by shifts of the module $V_0,$ and in the Grothendieck group 
  $$ [V_i] = [V_0] + q[V_0] + ... + q^{i} [V_0].$$ 
Therefore, the Grothendieck group is a cyclic $\Z[q,q^{-1}]$-module 
generated by $[V_0]=[\mathbf{1}].$ Moreover, $[H]=[V_{p-1}]=0,$ implying 
 $(1+q+\dots +q^{p-1})[\mathbf{1}]=0.$ The image of $\mathbf{1}$ is the 
unit element $1$ of the Grothendieck ring. Thus, the Grothendieck ring 
of $\mc{H}'$ is isomorphic to $R_p.$ 

$\square$ 

The category $\mc{H}'$ is symmetric triangulated. The shift functor 
$T$ acts on indecomposables as follows: 
 $$ T (V_i) \cong V_{p-2-i}\{i+1-p\} $$
Hence, $T^2(V) \cong V\{-p\}$ for any $V\in \mc{H}'.$ 

\vspace{0.1in} 

We just saw that objects of $\mc{H}'$ are easy to classify. In contrast, 
for all but first few primes $p$ classification of morphisms in     
$\mc{H}'$ is a wild problem. For instance, if $p>7$ the classification  
of morphisms from the object $W_0\otimes V_4$ to the object  
$$(W_1\otimes V_0) \oplus (W_2 \otimes V_2)\{-1\} \oplus  
 (W_3\otimes V_4)\{-2\}\oplus (W_4\otimes V_6)\{-3\}  
\oplus (W_5\otimes V_8)\{-4\},$$ 
modulo the action of the group $\prod_{i=0}^5\mathrm{GL}(W_i),$ 
where $W_i$ are finite-dimensional $k$-vector spaces, is equivalent 
to classification of representations of quiver with 6 vertices 
$v_i, 0 \le i \le 5$ and 5 arrows $v_i\to v_0, 1\le i\le 5,$ which is 
known to be wild. 

\vspace{0.2in}

{\bf Cyclic grading and category $\mc{H}.$}  
Instead of $\Z$-graded $H$-modules we could consider $\Z_m$-graded 
$H$-modules and the corresponding stable category, where $m\in \N.$ 
The prime $p$ should divide 
$m$ if we want the Grothendieck group of the category to be $R_p.$ 
The case of $\Z_p$-graded modules allows the greatest freedom 
(since a $\Z$-grading and a $\Z_{pr}$-grading collapse to a $\Z_p$-grading), 
and we denote by $\mc{H}$ the stable category of finite-dimensional 
$\Z_p$-graded $H$-modules. 

$\mc{H}$ has the following description which does not 
explicitly refer to a grading. Let $\widehat{H}$ be the Hopf algebra
over $k$ with primitive generators $X,Y$ and defining relations 
$$ [Y,X]=X, \hspace{0.2in} Y^p=Y, \hspace{0.2in} X^p=0.$$
 $\widehat{H}$ has dimension $p^2$ and is the restricted 
universal enveloping algebra of a p-Lie algebra with the 
same generators and relations. Consider the category  
$\widehat{H}\dmod$ of finite-dimensional left $\widehat{H}$-modules 
and its stable category $\widehat{H}\umod.$ Generator $Y$ acts 
semisimply on any finite-dimensional $H$-module $M,$ and we 
can decompose 
$$ M = \oplusop{r\in \mathbb{F}_p} M^r, \hspace{0.2in} 
   Ym = rm \hspace{0.1in} \mathrm{if} \hspace{0.1in} 
   m\in M^r.$$ 
Then $X M^r \subset M^{r+1}$ and $X^p=0.$ We see that
$\widehat{H}\dmod$ is equivalent to the category of 
$\Z_p$-graded $H$-modules, with the grading provided by 
the action of $Y.$ This equivalence descends to an equivalence 
between stable categories $\mc{H}$ and $\widehat{H}\umod$  
respecting their monoidal structures. 

Proposition~\ref{cyclo} holds for $\mc{H}$ as well, and other 
results can be restated for $\mc{H}$ in place of $\mc{H}',$ 
sometimes with only minor obvious adjustments.  
In particular, the Grothendieck ring of $\mc{H}$ is isomorphic to $R_p.$ 
One difference between $\mc{H}'$ and $\mc{H}$ is that the latter 
has only finitely many isomorphism classes of indecomposable modules: 
$V_i\{j\}$ over $0\le i\le p-2$ and $0\le j \le p-1.$  
Also, $T^2(V)\cong V$ for any $V$ in $\mc{H},$ since $V\{p\} =V.$ 

The above transformation of $\Z_p$-graded $H$-modules into 
modules over a larger Hopf algebra $\widehat{H}$ is a special 
case of Majid's bosonization construction [M]. 

\vspace{0.2in} 

{\bf Split Grothendieck rings of $\mc{H}'$ and $\mc{H}.$} 

We now look at how tensor products $V_i \otimes V_j$ break into 
direct sum of indecomposables in $\mc{H}'.$ It's convenient to 
slightly enlarge the category and allow $\frac{1}{2}\Z$-graded modules, 
for then we can balance $V_i$ about degree zero. We use the same 
notation $\mc{H}'$ for this bigger category. Define the balanced 
indecomposable modules by
  $$ \widetilde{V}_i =  V_i \{\frac{-i}{2}\}.$$ 

\begin{prop} 
The following isomorphism holds in the stable category: 
\begin{equation}\label{sumv}
 \widetilde{V}_i \otimes \widetilde{V}_j \cong \oplusop{m\in I} 
 \widetilde{V}_m, 
\end{equation}
where $I$ is the set $\{|i-j|, |i-j|+2, \dots, \mathrm{min}(i+j, 
   2p-i-j-4)\}$ and $0\le i,j\le p-2.$  
\end{prop} 

\emph{Proof:} The following isomorphism in $H\gmod$ (and, hence, 
in $\mc{H}'$) can be checked directly:  
$$\widetilde{V}_1\otimes \widetilde{V}_j\cong \widetilde{V}_{j-1} 
 \oplus \widetilde{V}_{j+1}, \hspace{0.2in} 1\le j \le p-2.$$
When $j=p-2$ the module $V_{j+1}$ is projective, so that in 
$\mc{H}'$ we have 
$$\widetilde{V}_1\otimes \widetilde{V}_{p-2} \cong \widetilde{V}_{p-3}.$$
This implies formula (\ref{sumv}) for $i=1.$ Since both 
$H\gmod$ and $\mc{H}'$ have the 
Krull-Schmidt property (unique direct sum decomposition), we can 
use induction on $i$ to check that (\ref{sumv}) holds whenever 
$i+j+2\le p.$ In this case the sum on the right extends to 
$i+j.$ 

\vspace{0.1in} 

Assume now $i\ge j$ and $i+j+2>p$ and apply the shift functor:  
$$T(\widetilde{V_i} \otimes \widetilde{V_j}) \cong 
(T\widetilde{V_i}) \otimes \widetilde{V_j}.$$
Since $T(\widetilde{V}_i) \cong \widetilde{V}_{p-2-i}\{\frac{-p}{2}\}$ 
and $(p-2-i)+j+2 \le p,$ we have 
 $$ (T\widetilde{V}_i) \otimes \widetilde{V}_j \cong 
   \oplusop{m\in J} \widetilde{V}_m\{\frac{-p}{2}\}, $$
where $J$ is the set $\{j+i+2-p, j+i+4-p, \dots, j+p-2-i\}.$ 
Applying $T^{-1}$ to both sides of this isomorphism nets us 
$$ \widetilde{V}_i \otimes \widetilde{V}_j \cong \oplusop{m\in I}\widetilde{V}_m $$
where $I=\{i-j,i-j+2,\dots, 2p-4-i-j\},$ which is the formula
(\ref{sumv}) in the case $i\ge j$ and $i+j+2>p.$ The remaining 
case $i<j$ follows by symmetry. 

$\square$ 

Define the split Grothendieck group $K_s(C)$ of an additive 
category $C$ as the abelian group with generators $[M],$ over 
all objects $M$ of $C,$ and relations $[M]=[M']+[M'']$ whenever 
$M\cong M'\oplus M''.$ If $C$ is additive monoidal, the split 
Grothendieck group is a ring. 

\begin{corollary} \label{corrv} 
The split Grothendieck ring of $\mc{H}'$ is naturally  
isomorphic to the level $p-2$ Verlinde algebra for $SU(2),$ 
where we consider the Verlinde algebra over the ground ring 
$\Z[q^{\pm \frac{1}{2}}].$ 
\end{corollary} 

Indeed, $K_s(\mc{H}')$ is a free 
$\Z[q^{\pm \frac{1}{2}}]$-module with the basis 
$[\widetilde{V}_0], \dots, [\widetilde{V}_{p-2}],$ 
and the structure coefficients of the multiplication are 
exactly those of the Verlinde algebra. $\square$  

\vspace{0.1in} 

Likewise, the split Grothendieck ring of $\mc{H}$ can 
be identified with the level $p-2$ Verlinde algebra 
for $SU(2),$ over the ground ring $\Z[q^{\frac{1}{2}}]/(q^n).$

\vspace{0.2in}

{\bf Characterizations of $H$ and $\mc{H}'.$}

{\it 1. Inseparable extensions and group schemes:} The Hopf algebra $H$ can be thought of as the 
algebra of functions on an infinitesimal group scheme. 
Suppose that $a\in k, \sqrt[p]{a}\notin k.$ 
The  field extension $k\subset K, K=k(y), y^p=a$ is inseparable and the 
Galois group of $K$ over $k$ is trivial. The does exist, however, a 
nontrivial infinitesimal Galois group (a finite group scheme) acting on 
$K.$ Namely,  $K$ is a comodule over the Hopf algebra $H,$ with the 
$k$-linear coaction 
$$\Delta_K : K \lra H\otimes_k K,   \hspace{0.2in} 
   \Delta_K(y^i) = \sum_{j=0}^i \binom{i}{j} X^{i-j}\otimes y^j.$$
This is analogous to the group action of $G$ on $K$ in case of an extension 
$k\subset K$ with a nontrivial Galois group $G$. The $k$-linear action 
$G\times K \lra K$ can be written as $k[G]\otimes_k K \lra K$ and, dually, 
as $K \lra k[G]^{\ast} \otimes_k K.$ Unlike the group algebra
$k[G],$ the Hopf algebra $H$ is self-dual, $H^{\ast}\cong H.$  

\vspace{0.1in} 

{\it 2. Restricted Lie algebra:} The Hopf algebra $H$ is the restricted universal 
enveloping algebra of the abelian $p$-Lie algebra $kX, X^{[p]}=0.$ 

\vspace{0.1in}

{\it 3.} The stable category $\mc{H}'$ is equivalent to each of the following 
categories: 

(I) The category whose objects are unbounded exact complexes of 
finite-dimensional graded projective $H$-modules and whose morphisms are 
homomorphisms of complexes modulo chain homotopies.

(II) The quotient of the bounded derived category $D^b(H\gmod)$ by the 
subcategory of complexes of projective modules. 

(III) The category of graded matrix factorizations over $k[x]$ with potential 
$x^p,$ modulo chain homotopies (although the obvious 
monoidal structure of $\mc{H}'$ is hard to describe in the 
language of matrix factorizations). 

Equivalences with (I) and (II) hold for any Frobenius algebra, 
see [Bu], [R]. Equivalence with (III) holds only for $H$ as above. 

\vspace{0.06in} 

Matrix factorizations are D-branes (boundary conditions) in 
Landau-Ginzburg models of the string theory. Hori, Iqbal and Vafa [HIV] 
extracted the structure coefficients of the multiplication 
in the SU(2) level $p-2$ Verlinde algebra (for all integer 
$p\ge 2$) from the category of matrix factorizations over 
$\C[x]$ with the potential $x^p.$ 
It's likely that Corollary~\ref{corrv} 
and its characteristic 0 analogue (Proposition~\ref{verch0} below) 
are directly related to [HIV].

\vspace{0.2in}
 
{\bf When $n$ is a power of $p$.} 
Our Hopf algebra $H$ has a natural generalization.  
 Let $n=p^m$ for some $m>0,$ and consider the Hopf algebra 
  $$H_{(m)}= k[X]/(X^n) , \hspace{0.2in}  
   \Delta(X)= X\otimes 1 + 1 \otimes X, \hspace{0.2in} 
  \mathrm{deg} X = 1,$$
 with the antipode $S(X)=-X$ and the counit 
 $\epsilon(X^i)= \delta_{i,0}.$  The stable category 
 $H_{(m)}\ugmod$ of finite-dimensional graded left 
 $H_{(m)}$-modules modulo projectives is symmetric tensor triangulated and its 
 Grothendieck group is naturally isomorphic to the 
 commutative ring 
  $$ R_n=\Z[q]/(1+q+\dots +q^{n-1}).$$ 
Some but not all of the results and observations above 
generalize to $H_{(m)}.$ The Grothendieck ring of $H_{(m)}\ugmod$ 
is $R_n.$ 
The Hopf algebra $H_{(m)}$ is not self-dual when $m>1,$ and 
the split Grothendieck ring of $H_{(m)}\ugmod$ has the multiplication 
table in the standard basis 
 $[\widetilde{V}_0], \dots, [\widetilde{V}_{n-2}]$ 
different from that of the Verlinde algebra.


\section{Nilpotent derivations} \label{done}
 
{\bf A conjecture.} 
Let's assume, optimistically, that we are on the right track and 
there does exists a categorification of quantum 3-manifold  
invariants with $\mc{H}$ or $\mc{H}'$ as the base category. 
We are already restricted to the case when $n$ is a prime $p.$ The 
normalized WRT invariant 
$\tau(M)$ of a closed 3-manifold $M$ lies in the ring $R_p,$ 
and its categorification $\Theta(M)$ should be an object of the 
category $\mc{H}$ and have at least the following properties: 

\begin{itemize} 
\item Multiplicativity: $\Theta(M_1 \# M_2) \cong 
\Theta(M_1) \otimes \Theta(M_2).$ 
\item Duality: $\Theta(\overline{M})\cong \Theta(M)^{\ast}.$
Notice that $\mc{H}$ has a contravariant duality functor 
$M\longmapsto M^{\ast}=\Hom_k(M,k).$ 
\item The invariant of a 3-sphere: $\Theta(S^3)\cong V_0,$ 
the unit object of $\mc{H}.$ 
\end{itemize} 

$\Theta(M)$ would split into direct sum of indecomposable 
modules $V_i\{j\}$ and their multiplicities $\tau^{i,j}(M)$ 
will be $\Z_+$-valued invariants of closed 3-manifolds, with 
$0\le i\le p-2$ and $1\le j \le p.$  
The invarint $\tau(M)$ will be
equal to the linear combination of $\tau^{i,j}(M)$: 
$$ \tau(M) = \sum_{i,j} [i]_q \hspace{0.03in} 
 q^j \tau^{i,j}(M), \hspace{0.3in} [i]_q=1+q+\dots +q^i.$$ 

The invariant of 
a 4-dimensional (decorated) cobordism $W$ between (decorated) 
3-manifolds $M_1, M_2$ should be a morphism 
$$\Theta(W): \Theta(M_1) \lra \Theta(M_2).$$ 
In particular, if the boundary of $W$ is empty, its invariant 
$\Theta(W)$ is an endomorphism of the unit object and thus 
an element of the field $k$ (possibly just a residue modulo $p$). 

\hspace{0.15in} 

When $\tau(M)$ is extended to 3-manifolds with boundaries, 
integrality cannot be completely retained. The space of 
conformal blocks $\tau(S)$ associated to an oriented closed 
surface $S$ has an almost integral structure, that of a module 
over the ring $R_{(p)}=\Z[\zeta^{\frac{1}{4}}, \frac{1}{p}]\subset\C,$ 
and the invariant of a 3-cobordism $M$ between $S_1$ and $S_2$ is an 
$R_{(p)}$-linear map $\tau(S_1) \lra \tau(S_2).$ 

We don't know how to categorify the ring $R_{(p)}.$ To get 
rid of powers of $\frac{1}{p},$ Gilmer [G] restricts the 
cobordisms category and only allows "targeted" cobordisms, 
where each component of a cobordism contains a component 
of the target surface. The Witten-Reshetikhin-Turaev functor, 
restricted to targeted cobordisms, can be renormalized to 
be integral, and the space of conformal blocks $\tau(S)$ 
associated to a surface becomes a free $\Z[\zeta^{\frac{1}{2}}]$ 
of $\Z[\zeta^{\frac{1}{4}}]$-module. 

Let us ignore the fractional powers of $\zeta$ and pretend that 
 $\tau(S)$ is a free $\Z[\zeta]$-module of the rank 
given by the Verlinde formula.  A categorification $\Theta(S)$ 
of $\tau(S)$ should be a triangulated category depending, up to 
equivalence, only on the genus $g$ of $S.$   
Moreover, $\Theta(S)$ should be a module-category over $\mc{H},$ 
with the Grothendieck group a free
$R_p$-module of rank equal to the dimension of the 
space of conformal blocks on a genus $g$ surface. 
Projective action of the mapping class group of $S$ should 
lift to a projective action on $\Theta(S).$ Since the surfaces might 
have to be decorated, one might have to restrict from the 
mapping class group to a subgroup of decoration-preserving 
automorphisms and to a projective action of that subgroup. 
The action on $\Theta(S)$ should be by exact and $\mc{H}$-linear 
functors. 

More generally, to a decorated 3-dimensional "targeted"
cobordism $M$ between surfaces $S_1$ and $S_2$ there should be assigned 
 an exact and $\mc{H}$-linear functor $\Theta(M): \Theta(S_1)\lra 
\Theta(S_2).$ To a decorated 4-dimensional cobordism $W$ with corners 
between 3D targeted cobordisms $M_1,M_2$ with 
$\partial_i(M_1)= \partial_i(M_2),$ $i=0,1,$
there should be assigned a natural transformation $\Theta(W)$ between 
functors $\Theta(M_1),\Theta(M_2).$ 

Categories $\Theta(S),$ functors $\Theta(M)$ and natural transformations 
$\Theta(W),$ over all $S$ and all allowed $M$ and $W,$ 
should assemble into a 2-functor from a 2-category of  
4-dimensional cobordisms with corners 
to the 2-category of triangulated categories. 
Objects of the latter are triangulated categories, 1-morphisms 
are exact functors and 2-morphisms are natural transformation. 

\vspace{0.2in} 

{\bf Module-categories over $\mc{H}.$}
To search for categories $\Theta(S)$ we need a large supply of triangulated 
module-categories over $\mc{H}.$ We explained in 
Section~\ref{hopfological} how to construct such categories from 
comodule-algebras over $H.$ A left $H$-comodule algebra is an 
associative unital algebra $A$ equipped with an algebra homomorphism 
$\Delta_A: A \lra H \otimes A$
which makes $A$ into an $H$-comodule. 
Since our $H$ is $\Z_p$-graded, we want $A$ to be $\Z_p$-graded as well, 
and $\Delta_A$ to be grading-preserving. We write 
$$\Delta_A(a) = \sum_{i=0}^{p-1} X^i \otimes \partial_i(a).$$
Comodule-algebra  axioms imply
$$\partial_0(a)=a, \hspace{0.1in} \partial_i = \frac{f_1^i}{i!}, 
\hspace{0.1in} 1\le i \le p-1, \hspace{0.1in} \mathrm{and} 
\hspace{0.1in} \partial_1^p=0,$$
and that $\partial_1$ is a derivation, 
$$\partial_1(ab) = \partial_1(a)b + a \partial_1(b).$$ 
Denote $\partial_1$ by $\partial.$ Since $\deg(\Delta_A)=0$ and 
$\deg(X)=1,$ derivation $\partial$ must have degree $-1.$ 

\vspace{0.1in} 

We see that a $\Z_p$-graded $H$-comodule algebra $A$ is an associative 
$\Z_p$-graded $k$-algebra with a degree $-1$ derivation $\partial: A\lra A,$
 $$ \partial(ab)= \partial(a)b + a\partial(b)$$
(notice the absence of $(-1)^{|a|}$ in the formula) such 
that $\partial^p=0.$ Clearly, $\partial(1)=0.$ 
The $H$-comodule structure of $A$ is given by
 $$ \Delta_A(a) = \sum_{i=0}^{p-1} X^i \otimes \frac{\partial^i(a)}{i!}= 
  \widetilde{\exp} (X\otimes\partial) a,$$
tilde denoting reduced exponential.

\vspace{0.1in} 

For any such $A,$ the category $A_H\ugmod$ of $\Z_p$-graded $A$-modules 
modulo homotopies is triangulated $\mc{H}$-module category. The 
Grothendieck group $K(A_H\ugmod)$ is naturally an $R_p$-module. 

For those $A$ that are smash products $H\# B,$ we can form 
the derived category $D(B,H)$ as well, by localizing quasi-isomorphisms 
in $A_H\ugmod.$ Again, $D(B,H)$ is a triangulated $\mc{H}$-module 
category, with its Grothendieck group naturally an $R_p$-module. 
Due to $H$ being self-dual (up to the grading reversal), to form the 
smash product $H\# B$ we just need $B$ to have a nilpotent 
derivation as above, but of degree $1$ rather than $-1.$  

There are plenty of such $A$'s and $B$'s, but we don't know
whether any of them could lead to a categorification of the conformal 
block spaces. Corollary~\ref{corrv} hints that a categorification 
of the conformal blocks on a torus comes from the Auslander 
algebra of $\mc{H}.$

\vspace{0.2in} 

{\bf Characteristic $0$ case.}
The category $\mc{H}$ admits a characteristic $0$ 
analogue, which we'll now describe. Pick $n>1$ and a
field  $k$ of characteristic $0$ which contains a primitive $n$-th 
root of unity $\zeta.$ Let $H_n$ be the $k$-algebra 
with generators $K,X$ and relations 
$$ K X = \zeta X K, \hspace{0.3in} K^n=1, \hspace{0.3in} 
X^n=0.$$
$H_n$ is a Hopf algebra (known as the Taft Hopf algebra) 
of dimension $n^2$ with the comultiplication 
$$ \Delta(K)= K \otimes K , \hspace{0.1in} \Delta(X)= X\otimes 1 + 
K \otimes X$$ 
and the antipode $S(K)= K^{-1},$  $S(X) = - K^{-1}X.$  
 
\vspace{0.1in} 

A left $H_n$-module $M$ restricts to a module over the subalgebra 
generated by $K,$ which is the group algebra of the finite cyclic 
group $\Z_n.$ Since the field $k$ contains a primitive $n$-th 
root of unity, we can decompose 
$$M=\oplusop{r\in \Z_n} M^{r}, \hspace{0.2in} Kx=\zeta^r x,
\hspace{0.1in} x\in M^r.$$ 
Then $X: M^r \to M^{r+1}$ is a nilpotent operator, $X^n=0,$ 
and $M$ is isomorphic to a direct sum of indecomposable 
modules $V_i\{j\}.$ The latter has a basis $b_j, b_{j+1}, 
\dots, b_{j+i}$ (the indices are taken modulo $n$) and 
the action 
$$ K b_r = \zeta^r b_r, \hspace{0.2in} Xb_r = b_{r+1}, 
\hspace{0.1in} j \le r < j+i, \hspace{0.2in}  Xb_{j+i}=0.$$ 
Shifting the grading of a module $M$ amounts to tensoring $M$ 
with one of the one-dimensional modules $V_0\{j\}.$ The 
module $V_0$ is the unit object $\mathbf{1}$ 
of the tensor category of $H_n$-modules. Modules $V_{n-1}\{j\}$ are
projective. Recall that $H_n\umod$ denotes the stable 
category of finite-dimensional $H_n$-modules. 

\begin{prop} The Grothendieck ring $K(H_n\umod)$ is 
isomorphic to the ring $R_n = \Z[q]/(1+q+\dots + q^{n-1}).$ 
\end{prop} 

To write down the split Grothendieck ring of $H_n\umod,$ we 
introduce fractional grading, allowing shifts by $\{\frac{1}{2}\}.$
We can interpret this by assuming that the ground field contains 
a square root of $\zeta$ and adding $K^{\pm \frac{1}{2}}$ to $H_n$ 
with the obvious relations $K^{\frac{1}{2}}X= 
\zeta^{\frac{1}{2}} XK^{\frac{1}{2}}$ etc.
Denote by $\widetilde{V}_i$ the module $V_i\{\frac{-i}{2}\}.$

\begin{prop} \label{verch0} 
The split Grothendieck ring of $H_n\umod$ 
is naturally isomorphic to the level $n-2$ Verlinde algebra for 
$SU(2).$ In $H_n\umod$ tensor products decompose as follows: 
\begin{equation*}
 \widetilde{V}_i \otimes \widetilde{V}_j \cong \oplusop{m\in I} 
 \widetilde{V}_m, 
\end{equation*}
where $I$ is the set $\{|i-j|, |i-j|+2, \dots, \mathrm{min}(i+j, 
 2n-i-j-4)\}$ and $0\le i,j\le n-2.$  
\end{prop} 

\begin{corollary} For any $Y_1,Y_2\in H_n\umod$ there is an 
isomorphism $Y_1\otimes Y_2\cong Y_2 \otimes Y_1.$ 
\end{corollary} 

This corollary by itself implies commutativity of the 
Grothendieck ring and the split Grothendieck ring of $H_n\umod.$ 
For $n>2,$ however, $H_n$ does not have a quasi-triangular 
structure, and the isomorphisms $Y_1\otimes Y_2 \cong Y_2\otimes Y_1$ 
cannot be made functorial.  

\vspace{0.1in} 

The category $H_n\umod$ seems to be the right characteristic 
zero analogue of the stable category $\mc{H}$ from 
Section~\ref{cyclotomic}. Our earlier conjecture about $\mc{H}$ 
being the base category for a categorification of quantum 
3-manifold invariants could be restated for $H_n\umod.$ 
One could look for triangulated categories $C$ over $H_n\umod$ 
with the Grothendieck group naturally isomorphic to the space 
of conformal blocks, etc. Unfortunately, the absence of 
functorial isomorphisms  $Y_1\otimes Y_2\cong Y_2 \otimes Y_1$ 
constitutes a serious foundational obstacle. Suppose we could 
assign an object $\Theta(M)$ to a 3-manifold $M$ such 
that, for instance, $\Theta(M\# N) \cong \Theta(M)\otimes \Theta(N)$ 
in a natural way. The permutation diffeomorphism 
$M\# N\cong N\# M$ would lead to isomorphisms  
$\Theta(M)\otimes \Theta(N)\cong \Theta(N)\otimes \Theta(M)$ 
which should be compatible with maps $\Theta(W): \Theta(M_1) \lra 
\Theta(M_2)$ assigned to 4-dimensional cobordisms $W$ between 
$M_1$ and $M_2.$ Nonexistence of functorial isomorphisms 
$Y_1\otimes Y_2\cong Y_2 \otimes Y_1$ apparently precludes 
such compatibility. 
Moreover, for composite $n$ the ring $R_n$ is not the 
cyclotomic ring $\Z[q]/\Psi_n(q).$ 

\vspace{0.1in} 

To end on a more uplifting note, we'll mention a 
relation between $H_n$ comodule-algebras and Kapranov's $q$-homological 
algebra. Kapranov [Ka] considered the category of $N$-complexes of objects of 
an abelian category $\mc{A},$ where an $N$-complex is a 
sequence of objects and morphisms with the composition 
of any $N$ consecutive morphisms being zero. If $\mc{A}$ is 
defined over a ground ring $k$ which contains a primitive 
$N$-th root of unity, an $N$-complex can be tensored with 
an $N$-complex of $k$-vector spaces. For such $\mc{A}$
 the notion of chain homotopy for $N$-complexes makes sense, 
so that the quotient category of complexes by homotopies can 
be formed. There exists a version of differential calculus 
with an $N$-complex taking place of the usual De Rham complex.   
Categories of $N$-complexes have been studied by  
Dubois-Violette [D-V1], [D-V2], Dubois-Violette and Kerner [D-VK], 
Kassel and Wambst [KW] and others. Bichon [B] pointed out a 
link between $N$-complexes and comodules over an infinite
version of the Taft Hopf algebra 
(with the relation $K^N=0$ omitted).  

\vspace{0.1in} 

Suppose that $\mc{A}$ is the category of modules over 
a $k$-algebra $B.$ The tensor product $A=H_N\otimes B$ has 
the trivial structure of an $H_N$-comodule algebra, with 
$H_N$ being the Taft Hopf algebra. The category 
of $A$-modules is (almost) the  
category of $N$-complexes over $B\dmod,$ the only 
difference being in the grading, which is $\Z_N$ in the 
former case and $\Z$ in the latter. Kapranov's 
homotopy category is essentially the homotopy category 
$A_{H_N}\umod$ and his homology groups $_pH_i$ can be 
inverpreted via indecomposables in the stable 
category $H_N\umod.$ In particular, the homotopy 
category of $N$-complexes of $B$-modules is triangulated 
and can be localized by inverting quasi-isomorphisms 
of $N$-complexes. 

\vspace{0.2in}

{\bf Acknowledgements:} 
I greatly benefited from discussions with Joseph Bernstein, Igor 
Frenkel, Greg Kuperberg and David Radford. This paper 
was written, for the most part, at the University of California in  
Davis. NSF provided partial support via grant DMS-0407784.

\vspace{0.2in} 

\section{References} 

\noindent 
[BCR] D.~J.~Benson, J.~F.~Carlson and J.~Rickard, 
Thick subcategories of the stable module category, 
\emph{Fund. Math.} {\bf 153} (1997), no. 1, 59--80.

\noindent
[BL] J.~Bernstein and V.~Lunts, Equivariant sheaves 
and functors, \emph{Lecture Notes in Math.} 1578, (1994) Springer-Verlag.  

\noindent 
[B] J.~Bichon, N-complexes et alg\'ebres de Hopf, 
\emph{C. R. Math. Acad. Sci. Paris} {\bf 337} (2003), no. 7, 
441--444.

\noindent 
[Bu] R.-O.~Buchweitz, Maximal Cohen-Macaulay modules and 
Tate-cohomology over Gorenstein rings, preprint (1986). 

\noindent 
[Ca]  J.~F.~Carlson, Modules and group algebras, 
Notes by Ruedi Suter. Lectures in Mathematics ETH Z\"urich, 
Birkhäuser Verlag, Basel, 1996.

\noindent 
[CF] L.~Crane and I.~B.~Frenkel, Four dimensional topological 
quantum field theory, Hopf categories, and the canonical bases, 
\emph{J. Math. Phys.} 35 (1994) 5136-5154, hep-th/9405183.

\noindent 
[DNR] S.~Dascalescu, C.~Nastasescu and S.~Raianu,  
Hopf algebras. An introduction. \emph{Pure and Applied 
Mathematics,} {\bf 235} Marcel Dekker Inc., New York-Basel, 2001. 

\noindent 
[D-V1] M.~Dubois-Violette, $d^N=0$ : Generalized homology, 
\emph{K-theory} {\bf 14} (1998) 371--404, arxiv q-alg/9710021. 

\noindent 
[D-V2] M.~Dubois-Violette, Lectures on differentials, generalized 
differentials and on some examples related to theoretical 
physics,  \emph{Quantum symmetries in theoretical physics 
and mathematics (Bariloche, 2000),}  59--94, Contemp. Math. 294, 
AMS 2002, arxiv math.QA/0005256. 

\noindent 
[D-VK] M.~Dubois-Violette and R.~Kerner, Universal q-differential 
calculus and q-analog of homological algebra, 
arxiv q-alg/9608026.

\noindent
[EGST] K.~Erdmann, E.~L.~Green, N.~Snashall and R.~Taillefer, 
Representation theory of the Drinfel'd doubles of a 
family of Hopf algebras, arxiv math.RT/0410017. 

\noindent 
[FS] R.~Farnsteiner and A.~Skowronski, 
Galois actions and blocks of tame infinitesimal group schemes, 
preprint 2004. 

\noindent 
[GM] S. I. Gelfand and Yu. I. Manin, Methods of homological algebra, 
2003. 

\noindent
[G] P.~Gilmer, Integrality for TQFTs, arxiv math.QA/0105059. 

\noindent 
[Ha] D.~Happel, Triangulated categories in the representation 
theory of finite dimensional algebras, London Math. Soc. Lect. 
Note Ser. 119, 1988, Cambridge University Press. 

\noindent 
[HIV] K.~Hori, A.~Iqbal and C.~Vafa, D-Branes and mirror symmetry, 
arxiv hep-th/0005247. 

\noindent 
[Ka] M.~Kapranov, On the q-analogue of homological algebra, 
arxiv q-alg/9611005. 

\noindent 
[KW] C.~Kassel and M.~Wambst, Alg\'ebre homologique des 
N-complexes et homologie de Hochschild aux racines de 
l'unit\'e, arxiv q-alg/9705001. 

\noindent 
[KN] B.~Keller and A.~Neeman, The connection between May's axioms 
for a triangulated tensor product and Happel's description of the derived 
category of the quiver $D_4,$ \emph{Documenta Mathematica} {\bf 7} (2002), 535-560. 

\noindent 
[Kh] M.~Khovanov, A categorification of the Jones polynomial, 
\emph{Duke Math. J.,} {\bf 101}  (2000),  no. 3, 359--426,
arxiv math.QA/9908171. 

\noindent 
[Ku] G.~Kuperberg, Non-involutory Hopf algebras and 3-manifold 
invariants, \emph{Duke Math J.,} {\bf 84} (1), 83-129 
(1996), arXiv q-alg/9712047. 

\noindent 
[LS] R.~G.~Larson and M.~E.~Sweedler, An associative orthogonal 
bilinear form for Hopf algebras, \emph{Amer. J. Math.} {\bf 91} 
(1969), 75--94.

\noindent
[Le] T.~T.~G.~Le, Strong integrality of quantum invariants of 
3-manifolds, math.GT/0512433.  

\noindent 
[M] S.~Majid, Comments on bosonisation and biproduct, 
arxiv q-alg/9512028. 

\noindent 
[MR] G.~Masbaum and J.~Roberts, A simple proof of integrality of 
quantum invariants at prime roots of unity, \emph{Math. Proc. Cam. Phil. Soc.}
 121 (1997) 443-454. 

\noindent 
[Ma] J.~P.~May, The additivity of traces in triangulated 
categories, \emph{Advances in Math.} {\bf 163}, 34--73 (2001). 

\noindent 
[MVW] C.~Mazza, V.~Voevodsky and C.~Weibel, Notes on 
motivic cohomology, available at 
http://math.rutgers.edu/$\tilde{\quad}$weibel  

\noindent 
[Mo] S.~Montgomery, Hopf algebras and their actions on rings, 
AMS 1993. 

\noindent 
[Mu] H.~Murakami, Quantum SO(3) invariants dominate the SU(2) 
invariant of Casson and Walker, \emph{Math. Proc. Cam. Phil. Soc.}
 117 (1995) 237--249. 

\noindent 
[Ra] D.~E.~Radford, The trace function and Hopf algebras, 
\emph{Journal of algebra} {\bf 163} 583--622, (1994). 

\noindent 
[RT] N.~Reshetikhin and V.~Turaev, Invariants of 3-manifolds via 
link polynomials and quantum groups, \emph{Invent. Math.} 
{\bf 103}, 547--597 (1991). 

\noindent 
[R] J.~Rickard, Derived categories and stable equivalence, 
\emph{Journal of Pure Appl. Alg.} {\bf 61} (1989), 303-317. 

\noindent 
[St] N.~P.~Strickland, Axiomatic stable homotopy -- a survey, 
preprint arXiv math.AT/0307143. 

\noindent 
[V] V.~Voevodsky, Homology of schemes and covariant 
motives, PhD thesis, Harvard University, 1992. 

\noindent
[We] C.~Weibel, An introduction to homological algebra, 1996. 

\noindent 
[W] E.~Witten, Quantum field theory and the Jones polynomial, 
\emph{Comm. Math. Phys.} {\bf 121}, 351--399 (1989). 

\vspace{0.3in} 

\noindent
Department of Mathematics, Columbia University, New York, NY 10027. 

\noindent
khovanov@math.columbia.edu

\end{document}